\documentclass[12pt]{amsart}
\usepackage{amssymb}

\let\frak\mathfrak
\let\Bbb\mathbb

\def\>{\relax\ifmmode\mskip.666667\thinmuskip\relax\else\kern.111111em\fi}
\def\<{\relax\ifmmode\mskip-.333333\thinmuskip\relax\else\kern-.0555556em\fi}
\def\vsk#1>{\vskip#1\baselineskip}
\def\vv#1>{\vadjust{\vsk#1>}\ignorespaces}
\def\vvn#1>{\vadjust{\nobreak\vsk#1>\nobreak}\ignorespaces}

\let\Medskip\medskip
\def\medskip{\par\Medskip}
\let\Bigskip\bigskip
\def\bigskip{\par\Bigskip}

\let\Maketitle\maketitle
\def\maketitle{\Maketitle\thispagestyle{empty}\let\maketitle\empty}

\newtheorem{thm}{Theorem}[section]
\newtheorem{cor}[thm]{Corollary}

\numberwithin{equation}{section}

\theoremstyle{definition}

\let\mc\mathcal
\let\nc\newcommand

\nc{\on}{\operatorname}
\nc{\Z}{{\mathbb Z}}
\nc{\C}{{\mathbb C}}
\nc{\N}{{\mathbb N}}
\nc{\pone}{{\mathbb C}{\mathbb P}^1}
\nc{\arr}{\rightarrow}
\nc{\larr}{\longrightarrow}
\nc{\al}{\alpha}
\nc{\W}{{\mc W}}
\nc{\la}{\lambda}
\nc{\su}{\widehat{{\mathfrak sl}}_2}
\nc{\g}{{\mathfrak g}}
\nc{\h}{{\mathfrak h}}
\nc{\m}{{\mathfrak m}}
\nc{\n}{{\mathfrak n}}
\nc{\Gm}{\Gamma}
\nc{\La}{\Lambda}
\nc{\gl}{\widehat{\mathfrak{gl}_2}}
\nc{\bi}{\bibitem}
\nc{\om}{\omega}
\nc{\Res}{\on{Res}}
\nc{\gm}{\gamma}
\nc{\Om}{\Omega}

\def\Res{\on{Res}}

\def\Wr{\on{Wr}}

\def\B{{\mc B}}

\def\M{{\mc M}}
\def\O{{\mc O}}

\let\leq\leqslant

\nc{\gln}{\mathfrak{gl}_N}
\nc{\sln}{\mathfrak{sl}_N}

\nc{\glN}{\mathfrak{gl}_N}
\nc{\slN}{\mathfrak{sl}_N}

\def\beq{\begin{equation}}
\def\eeq{\end{equation}}
\def\be{\begin{equation*}}
\def\ee{\end{equation*}}

\nc{\bean}{\begin{eqnarray}}
\nc{\eean}{\end{eqnarray}}
\nc{\bea}{\begin{eqnarray*}}
\nc{\eea}{\end{eqnarray*}}
\nc{\bs}{\boldsymbol}
\nc{\Ref}[1]{{\rm(\ref{#1})}}

\nc{\Wl}{\Wr_{\bs\la}}
\nc{\Ol}{\Om_{\bs\la}}
\nc{\Mla}{\M_{\bs\la,\bs a}}
\nc{\p}{\partial}
\nc{\Bl}{\B_{\bs\la}}
\nc{\Bla}{\B_{\bs\la,\bs a}}
\nc{\Ola}{\O_{\bs\la,\bs a}}

\nc{\OL}{\mc O_{\bs\La,\bs\la,\bs b}}

\nc{\Bm}{\mc B_{m^{\bs\La}_{\bs\la}}}
\nc{\R}{\Bbb R}

\nc{\s}{{\rm sing}}
\nc{\Oll}{{\Omega_{\bs\la}}}

\nc{\Tee}{\mathcal S}
\nc{\un}{U(\n_-)}

\nc{\ep}{\epsilon}
\nc{\slt}{\frak{sl}_2}
\nc{\sltu}{\frak{sl}_2[t]}
\nc{\Pk}{{(P^1)^k}}
\nc{\PN}{{(P^1)^N}}
\nc{\Ck}{{\mathcal C_k(z)}}
\nc{\CN}{\mathcal C_N(z)}
\nc{\AN}{{\mathcal A^N}}
\nc{\Ak}{{\mathcal A^k}}
\nc{\mP}{{\tilde{\mathcal P}}}
\nc{\mD}{{\tilde{\mathcal D}}}
\nc{\zz}{{\bs z}}
\nc{\yy}{{\bs y}}
\nc{\TT}{{\bs t}}

\begin{document}

\title
[A Selberg Integral
Type Formula]
{A Selberg Integral
Type Formula for an $\slt$ One-Dimensional Space of Conformal Blocks}

\author[A.\,Varchenko]
{A.\,Varchenko$\>^{1}$}

\thanks{${}^1$\ Supported in part by NSF grant DMS-0555327}

\maketitle

\begin{center}

\medskip
{\it Department of Mathematics, University of North Carolina
at Chapel Hill\\ Chapel Hill, NC 27599-3250, USA\/}
\end{center}

\medskip
\begin{abstract}
For  distinct complex numbers 
$z_1,\dots,z_{2N}$, we
give  a polynomial $P(y_1,\dots,y_{2N})$ in the variables
$y_1,\dots,y_{2N}$, which is homogeneous of degree $N$,
linear with respect to each variable,  
$\slt$-invariant
with respect to a natural $\slt$-action, and is of order $N-1$ at
$(y_1,\dots,y_{2N})=(z_1,\dots,z_{2N})$.

We give also a Selberg integral type formula for the associated
one-dimensional space of conformal blocks.

\medskip
\noindent
2000 Math. Subj. Class. Primary 81T40, 33C70; Secondary 32S40, 52B30.

\medskip
\noindent
Key words and phrases. Conformal blocks, invariant polynomials.
\end{abstract}

\maketitle
\section{Introduction}

According to a general principle in \cite{MV}, if a KZ-type equation
has a one-dimensional space of solutions, then the hypergeometric integrals
representing the solutions can be calculated explicitly, see demonstrations
of that principle in \cite{TV}, \cite{FSV}, \cite{W}. In this note we give another
example of that type. 

We consider the bundle of the $\slt$ conformal blocks at level one on
the Riemann sphere. That bundle is of rank one. 
 Our first result is a formula for a generator of a fiber
of that bundle, see Theorem \ref{thm formula}.  Namely, for
distinct complex numbers $\zz=(z_1,\dots,z_{2N})$, we  
give a polynomial $P(y_1,\dots,y_{2N})$ in the variables
$\yy=(y_1,\dots,y_{2N})$, which is homogeneous of degree $N$,
 linear with respect to each variable,  $\slt$-invariant
with respect to a natural $\slt$-action, and is of order $N-1$ at
$\yy=\zz$.

The conformal block bundle has a KZ connection.
The flat sections of the KZ connection  have
representations in terms of multidimensional hypergeometric integrals,
see \cite{SV}, \cite{FSV1}, \cite{FSV2}. The formula of Theorem 
\ref{thm formula} for a generator of
a conformal block space allows us to calculate 
those integrals
explicitly, see Theorem \ref{thm main}.
 The result is a Selberg integral type formula \Ref{selberg}.

\medskip
The author thanks I. Dolgachev, 
whose interest to this subject had triggered
the writing this note, and  R. Rimanyi, V. Schechtman for useful discussions.

\section{ Spaces of Conformal Blocks}
\subsection{Conformal blocks}
\label{subsec Conf Bl}

Consider the complex Lie algebra $\slt$ with
generators $e,f,h$ and relations $[e,f]=h, [h,e]=2e, [h,f]=-2f$. For a
nonnegative integer $m$, denote by $V_m$ the irreducible $\slt$-module
with highest weight $m$. 

\medskip

 Let $\ell, m_1,\dots,m_n, N$ be given nonnegative integers such that
\bean
\label{integer cond}
\ell > 0\ ,
\qquad
0\ \leq m_1, \dots, m_n,\ m_1+\dots+m_n-2N \ \leq \ell \ .
\eean
Denote \ $p =  \ell + 1+ 2N - m_1 - \dots - m_n$,
\bea
\label{L product}
V\ = \ V_{m_1}\otimes \dots \otimes V_{m_n}\ .
\eea
For each $a=1,\dots,n$, denote by
$e_a:V\to V$ the linear operator  acting as
 $e$ on the $a$-th factor and as the identity on all the other factors.

Let $\zz=(z_1,\dots,z_n)$ be a collection of distinct complex numbers.
Denote
\bean
\label{def of conf blocks}
W(\zz)\ = \ \{ v \in V\ | \ hv=(\sum_{a=1}^nm_a - 2N)v,\ ev=0,\ 
(\sum_{a=1}^nz_ae_a)^pv=0\}\ .
\eean
The vector
space $W(\zz)$ 
is called the {\it space of conformal blocks at level 
$\ell$}, see \cite{FSV1}, \cite{FSV2}.

\medskip
\noindent
{\bf Remark.}
This definition is nonstandard. Usually the space of conformal
 blocks  is defined if one has
$n$ distinct points on a Riemann surface 
and $n$ irreducible  representations of an affine
 Lie algebra, see \cite{KL}. If the Riemann surface
 is the Riemann sphere,
 then one can describe the space of conformal blocks in terms
 of finite dimensional representations of the corresponding finite
dimensional Lie algebra. That
 description is one of two main results of \cite{FSV1} and \cite{FSV2}. 
We take that description as our definition.

\subsection{KZ connection}
Denote
\bea
X_n\ = \ \{ \zz=(z_1,\dots,z_n) \in 
\C^n\ | \ z_a\neq z_b\ {\rm for\ all }\ a\neq b\}\ .
\eea
The trivial vector  bundle 
$\eta  : V\times X_n\to X_n$
has a KZ connection,
\bea
\frac{\partial}{\partial z_a} \ - \ \frac 1{\ell+2}\,
\sum_{b\neq a} \frac {\Omega^{(a,b)}}{z_a-z_b}\ ,
\qquad
a=1,\dots,n\ ,
\eea
where $\Omega = \frac 12\,h\otimes h + e\otimes f + f\otimes e$\
and $\Omega^{(a,b)}:V\to V$ is the linear operator acting as $\Omega $ on the 
$a$-th and $b$-th factors and as the identity on all the other factors.

Consider the subbundle of conformal blocks
with fiber $W(\zz)\subset V$. It is well known that
this subbundle is invariant with respect to the KZ connection
\cite{KZ}.

\section{Conformal blocks at level one}
In the rest of the paper we assume that
\bean
\label{number assumption}
\ell\ =\ 1\ ,
\qquad
 n\ =\ 2N\ ,
\qquad
 m_1 = \dots = m_{2N} = 1\ .
\eean
These numbers satisfy conditions \Ref{integer cond}.
Then
\bea
W(z_1,\dots,z_{2N})\ =\ \{ v \in (V_1)^{\otimes 2N}\ |\
hv=0,\ 
ev=0,\ (\sum_{a=1}^nz_ae_a)^2v=0\}\ .
\eea

\begin{thm}
\label{thm on dim}
Under assumptions \Ref{number assumption},
 the space of conformal blocks is one-dimensional.
\end{thm}

\begin{proof}
The fusion ring  of $\slt$ at level $\ell=1$ is a free
$\Z$-module with generators $[V_0], [V_1]$ and commutative associative multiplication:\
\bea
[V_0]\cdot [V_0] = [V_0],
\qquad 
[V_0]\cdot [V_1] = [V_1],
\qquad
[V_1]\cdot [V_1] = [V_0].
\eea 
The dimension of
$W(z_1,\dots,z_{2N})$ is the coefficients of $[V_0]$ in the decomposition of
$[V_1]^{2N}$ in terms of the generators. Clearly $[V_1]^{2N} = [V_0]$.
\end{proof}

\medskip

Let us 
realize the tensor product $(V_1)^{\otimes 2N}$ as the vector space
of polynomials
\bea
p(\yy)\ = \ p(y_1,\dots,y_{2N})
\eea 
of degree not greater than one
with respect to each variable $y_1,\dots,y_{2N}$. 
The Lie algebra $\slt$  acts on this
space in the standard way, in particular, 
$e$ acts as $\sum_{a=1}^{2N} \partial/
\partial y_a$.

By \cite{R}, Theorem 4.3  (cf. \cite{LV}, Lemma 1.3), under assumptions
\Ref{number assumption},
the subspace  $W(\zz) \subset (V_1)^{\otimes 2N}$
consists of polynomials  $p(\yy)$,
which are $\slt$-invariant, homogeneous of degree $N$
and of order at least $N-1$ at $\yy=\zz$.
(The formulation of
Theorem 4.3 in \cite{R} has a misprint: $\phi$ should vanish to order $J-k$
rather than $J - k - 1$).

Introduce a polynomial $P(\yy;\zz)$ in the variables
$\yy$ depending on the parameters $\zz$,
\bean
\label{formula main}
P(\yy;\zz)\ = \ 
\det_{1\leq a,b\leq N} \left(\frac{y_{a}-y_{N+b}}{z_{a}-z_{N+b}}\right)\ =\
\sum_{\sigma\in S_N}\,(-1)^\sigma\prod_{a=1}^N \,
\frac{y_{\sigma_a}-y_{N+a}}{z_{\sigma_a}-z_{N+a}}\ .
\eean

\begin{thm}
\label{thm formula}
For a fixed $\zz$, the polynomial $P(\yy;\zz)$ in the variables $\yy$
is $\slt$-invariant, is
homogeneous of degree $N$, has degree not greater than
one with respect to each variable $y_1,\dots,y_{2N}$ and has order $N-1$ at
$\yy=\zz$; therefore, this polynomial
is a generator of the space  $W(\zz)$ of conformal blocks.
\end{thm}

\begin{proof} 
Each difference $y_i-y_j$ is $\slt$-invariant. Hence $P(\yy,\zz)$ is $\slt$-invariant.
We have
$P(\zz;\zz)=0$ since it is the determinant of a matrix with all entries
equal to one. The fact that $P(\yy;\zz)$ has zero of order $N-1$ at $\yy=\zz$ 
 is proved similarly.
\end{proof}

Denote
\bean
\label{A}
\phantom{aaa}
A(\zz)\ =\
\prod_{1\leq a\leq N < b\leq 2N} (z_b-z_a)^{1/2} \
\prod_{1\leq a < b\leq N} (z_b-z_a)^{-1/2}\ 
\prod_{N <a < b\leq 2N} (z_b-z_a)^{-1/2}\ .
\eean

\begin{thm}
\label{thm section s}

Consider the trivial bundle $\eta : (V_1)^{\otimes 2N} \times X_{2N}\to X_{2N}$
and its section 
\bean
\label{section}
s\ :\ \zz \mapsto \ A(\zz)P(\yy;\zz) \in W(\zz)\ .
\eean
 Then this is a flat section of the KZ connection for $\ell=1$.
\end{thm}

The theorem can be proved by a direct calculation. A different proof see in
Section \ref{sec proofs}.

\medskip
\noindent
{\bf Remark.} The flat
section $s$ is multivalued with the monodromy equal to
-1 around each hyperplane
$z_a=z_b$. Hence, the KZ connection on the  $\slt$ conformal block bundle
at level one is unitarizable. The unitarity of the KZ connection on the 
$\slt$ conformal block bundle at any level is proved in \cite{R}, cf \cite{LV}.

\section{An integral representation for conformal blocks}
\label{sec int repn}

\subsection{The master and weight functions}
\label{sec master and weight}
Introduce a scalar function
\bea
\Phi(\TT;\zz) &=& 
\Phi(t_1,\dots,t_N;z_1,\dots,z_{2N})\ 
\\
&=& 
\prod_{1\leq a<b\leq 2N}
(z_a-z_b)^{1/6}\
\prod_{1\leq i<j\leq N}
(t_j-t_i)^{2/3}\
\prod_{i=1}^N\prod_{a=1}^{2N}\
(t_i-z_a)^{-1/3} \ 
\eea
and a $(V_1)^{\otimes 2N}[0]$-valued rational function,
\bea
\omega(\TT;\zz)\ = \ \sum_{1\leq a_1<\dots <a_N\leq 2N}\
\sum_{\sigma\in S_N}\ \prod_{i=1}^N\ \frac {y_{a_i}}{t_{\sigma_i}-z_{a_i}}\ .
\eea
The functions $\Phi$ and $\omega$ are called the {\it master} and {\it
weight} functions, respectively, see \cite{SV}.

\subsection{The local system}
\label{sec loc system}
Denote
\bea
Y_{2N}
&=&
  \{ (\TT;\zz)=
(t_1,\dots,t_N;z_1,\dots,z_{2N}) \in \C^{3N}\ | 
\\
&&
\phantom{aaaaaaa}
\ z_a\neq z_b\ 
{\rm for\ all }\ a< b;\ 
 t_i\neq t_j\ 
 {\rm for\ all }\ i< j;\ 
z_a\neq t_i\ 
{\rm for\ all }\ a, i\}\ .
\eea
The master function $\Phi$ defines on $Y_{2N}$
a one-dimensional local system $\mc
L$. The horizontal sections of $\mc L$
are generated by the univalued branches of
the multivalued holomorphic function $\Phi$.

The projection $\tau : Y_{2N} \to X_{2N}, \, (\TT;\zz) \mapsto \zz$ is
topologically trivial.  Let $\tau_N$ be the associated {\it
homological} vector
bundle with fiber $H_N (\tau^{-1}(\zz), \mc L|_{\tau^{-1}(\zz)})$.
  The vector bundle $\tau_N$ has a canonical
Gauss-Manin connection.

\subsection{The integral representation}
\label{Sec int repn}

\begin{thm}[\cite{SV},\cite{FSV1},\cite{FSV2}]
\label{thm SV}
Let $\gamma$ be a horizontal section of the homological bundle $\tau_N$.
Then 
\bean
\label{Integral}
I_{\gamma}(\zz)\ = \ 
\int_{\gamma(\zz)}\,\Phi(\TT;\zz)\,\omega(\TT;\zz)\,
dt_1\wedge\dots\wedge dt_N\ 
\eean
is a horizontal section of the KZ connection. Moreover, this section
takes values in the
conformal block spaces.
\end{thm}

\begin{cor}
For any horizontal section $\gamma$, there exists $c_\gamma\in \C$ such that
$I_\gamma(\zz) = c_\gamma s(\zz)$, where  $s$ is the horizontal
section defined in
\Ref{section}.
\end{cor}

\section{An example of a horizontal family $\gamma$}
\label{sec ex of horiz family}

\subsection {Euler's beta function} 
\label{Euler beta function} 
Let $z_a< z_b$ be real numbers and $\alpha,\beta$ positive
numbers. For $t\in (z_a,z_b)$, we fix ${\rm arg} (t-z_a)=0$,\ ${\rm arg}
(t-z_b)=\pi$, ${\rm arg} (z_b-z_a)=0$.  Then
\bean
\label{beta formula}
\int_{z_a}^{z_b}\,(t-z_a)^{\alpha-1}(t-z_b)^{\beta-1}\,dt
\ =\ - e^{\pi i \beta}\,\frac { \Gamma(\al) \Gamma(\beta)}
{\Gamma(\alpha+\beta)}\, (z_b-z_a)^{\alpha+\beta-1}\ .
\eean
The right hand side of 
\Ref{beta formula} 
is a holomorphic function of
$\alpha,\beta \in \C - \{0,-1,\dots\}$. We define the integral in
the left hand side of
\Ref{beta formula} for $\alpha,\beta \in \C - \{0,-1,\dots\}$ by
analytic continuation.

\medskip
Fix $\alpha,\beta \in \C$. The function
$(t-z_a)^{\alpha-1}(t-z_b)^{\beta-1}$ defines on $\C-\{z_a,z_b\}$ a one-dimensional
local system $\ell$, whose sections are generated by the univalued
branches of that function. It is easy to see that if $\alpha,\beta,
\alpha+\beta \in \C-\Z$,
then there is a unique
cycle $\gamma(z_a,z_b;\alpha,\beta)
\in H_1(\C-\{z_a,z_b\}, \ell)$ such that
\bea
\int_{\gamma(z_a,z_b;\alpha,\beta)}\,(t-z_a)^{\alpha-1}(t-z_b)^{\beta-1}\,dt\
=\
\int_{z_a}^{z_b}\,(t-z_a)^{\alpha-1}(t-z_b)^{\beta-1}\,dt\ .
\eea
This cycle will be called a {\it Pochhammer cycle}.

\subsection{An example of a horizontal section  $\gamma$ of the
homological bundle $\tau_N$} 
\label{ex horiz}

Assume that $\zz \in \R^{2N}$ and 
\bean
\label{real z}
z_1< z_{N+1}
< z_2 < z_{N+2}<\dots < z_N < z_{2N}\ .
\eean
First, we define $\gamma(\zz)$ as an oriented product of intervals, 
\bean
\label{cycle}
\gamma(\zz) \ = \ \{ 
(\TT;\zz)\ | \
t_a \in (z_a,z_{N+a}), \ a=1,\dots,N\}\ 
\eean
with the standard orientation of each of the intervals $(z_a,z_{N+a})$

To define the integral in \Ref{Integral}, we fix on $\gamma(\zz)$ the
arguments of all factors of the master function $\Phi$ as follows.  We
set ${\rm arg}(t_j-t_i)=0$ for $j>i$.  We set ${\rm arg}(z_a-z_b)=0$
if $z_a>z_b$ and ${\rm arg}(z_a-z_b)=\pi$ if $z_a<z_b$.  We set ${\rm
arg}(t_i-z_a)=0$ if $t_i>z_a$ and ${\rm arg}(t_i-z_a)=\pi$ if
$t_i<z_a$.  This assignment of arguments determines the integral
$I_\gamma(\zz)$ in \Ref{Integral} and determines a horizontal section
$\gamma$ of the homological bundle $\tau_N$ for real $\zz$ satisfying
the above conditions. We extend this horizontal section to other values of 
$\zz$ by continuity. 
 
\medskip
\noindent
{\bf Remark.}
Strictly
speaking we have defined $\gamma(\zz)$ as a cell with 
a coefficient in $\mc
L$. The boundary of that cell lies in the union of hyperplanes
of the singularities of the master function $\Phi$. Nevertheless,
using the
remark in Section \ref{Euler beta function}, we can represent the same
function  $I_\gamma(\zz)$ as an integral over the product of the
corresponding Pochhammer cycles 
\bea
\gamma(z_1,z_{N+1};-1/3,-1/3)\times \dots \times
\gamma(z_N,z_{2N};-1/3,-1/3)\ ,
\eea
 that is, we can represent $I_\gamma(\zz)$ 
as an integral over an element
of $H_N (\tau^{-1}(\zz), \mc L|_{\tau^{-1}(\zz)})$ and that element
depends on $\zz$ horizontally.

\subsection{A Selberg integral type formula}
 Assume that $\zz \in \R^N$  satisfies  \Ref{real z}.
Let $s : \zz \mapsto A(\zz)P(\yy;\zz)$ 
be the section defined 
in \Ref{section}. The function $A(\zz)$ is multivalued.
We fix its univalued branch over the set of  $\zz$'s satisfying
\Ref{real z} by setting
${\rm arg}(z_a-z_b)=0$ if $z_a>z_b$ and
${\rm arg}(z_a-z_b)=\pi$ if $z_a<z_b$.

\begin{thm}
\label{thm main}
 Let $\zz \in \R^N$ satisfy \Ref{real z}.
Let $s$ 
be the section  defined 
in \Ref{section}.
Let $\gamma$ be the section of  $\tau_N$ defined in Section
\ref{ex horiz}. 
Then\ $I_\gamma(\zz)\, =\, C\, s(\zz)$, where
\bean
\label{constant}
C\ = \  e^{-\pi i N^2/3}\
(3 e^{\pi i/6}  \Gamma(2/3)^3 \sin(\pi/3) /\pi)^N\ .
\eean
\end{thm}

In particular, comparing the coefficients of the monomial
$y_1\dots y_N$ in the
right and left hand sides of the equation $I_\gamma(\zz)\, =\, C\, s(\zz)$,
we get the following formula,
\bean
\label{selberg}
&&{}
\\
&&
\int_{z_1}^{z_{N+1}}\!\!\!\!\! 
\dots \!
\int_{z_N}^{z_{2N}}\!\!\!
\prod_{1\leq i<j\leq N}
(t_j-t_i)^{2/3}
\prod_{i=1}^N\prod_{a=1}^{2N}\
(t_i-z_a)^{-1/3} \times
\phantom{aaa}
\notag
\\
&&
\phantom{aa}
\times \sum_{\sigma\in S_N} \prod_{a=1}^N\ \frac {1}{t_{\sigma_a}-z_{a}}\ 
dt_1\wedge\dots\wedge dt_N\ 
\ = \  C\!\!
\prod_{1\leq a\leq N < b\leq 2N} (z_b-z_a)^{1/3} \times
\notag
\\
&&
\phantom{aaaaa}
\times
\prod_{1\leq a < b\leq N} (z_b-z_a)^{-2/3}\!\!
\prod_{N <a < b\leq 2N} (z_b-z_a)^{-2/3}
\sum_{\sigma\in S_N}\,(-1)^\sigma\prod_{a=1}^N \,
\frac{1}{z_{\sigma_a}-z_{N+a}}\ .
\notag
\eean
This is a Selberg integral type formula.

\subsection{A proof of Theorems \ref{thm section s} and \ref{thm main}}
\label{sec proofs}
The KZ connection on the bundle $\eta$
has regular singularities. Therefore, any horizontal section $\tilde s$
of the conformal block subbundle has the form
$\tilde s : \zz \mapsto \tilde A(\zz) P(\yy;\zz)$, where
\bea
\tilde A(\zz)\ =\
\prod_{1\leq a < b\leq 2N} (z_b-z_a)^{\alpha_{a,b}} 
\eea
for suitable numbers $\alpha_{a,b}$. To prove Theorem \ref{thm section s}
we need to show that the numbers $\alpha_{a,b}$ are given by formula
\Ref{A}.

Assume that $z_{N+1}-z_1,\dots,z_{2N}-z_N$ all tend to zero. Then
\bea
&&
I_\gamma(\zz)\ = \ e^{-\pi i\frac{N(N-1)}{12}}\ \times
\phantom{aaaaaaaaaaaaaaa}
\\
&&
\phantom{aaaa}
\prod_{a=1}^N \,(z_a-z_{N+a})^{1/6}
\int_{z_a}^{z_{N+a}}\! (t_a-z_a)^{-1/3}
(t_a-z_{N+a})^{-1/3} 
(\frac{y_a}{t-z_a} + \frac{y_{N+a}}{t-z_{N+a}}) \,dt_a + \dots\ .
\eea
where the dots denote the higher order terms. Calculating the integrals we get
\bean
\label{integral}
&&{}
\\
\phantom{aa}
I_\gamma(\zz) &=&  e^{-\pi i\frac{N(N-1)}{12}}e^{\pi i \frac{5N}6}
(-1)^N \frac {\Gamma(-1/3)^N \Gamma(2/3)^N}{\Gamma(1/3)^N}  
\prod_{a=1}^N \ (z_{N+a}-z_a)^{-1/2} (y_a-y_{N+a}) + \dots
\notag
\\
&=&
e^{-\pi i\frac{N(N-1)}{12}}e^{\pi i \frac{5N}6}
(3\Gamma(2/3)^{3}\sin(\pi/3)/\pi)^N\
\prod_{a=1}^N \ (z_{N+a}-z_a)^{-1/2} (y_a-y_{N+a}) + \dots\ .
\notag
\eean
Comparing these asymptotics with the asymptotics of $\tilde A(\zz)P(\yy;\zz)$
we conclude that 
\linebreak
$\alpha_{a,N+a}= 1/2$ for all $a$. This statement is 
in agreement with formula \Ref{A}. Since the order of vanishing of
conformal blocks as $z_a-z_b$ tends to zero is the same for all pairs
$(a,b)$, we concllude that  $\tilde A(\zz)= A(\zz)$. Theorem \ref{thm section s}
is proved.

To prove Theorem \ref{thm main} we need to calculate the asymptotics of
$A(\zz) P(\yy;\zz)$ as 
\linebreak
$z_{N+1}-z_1,\dots,z_{2N}-z_N$ all
tend to zero and compare them
 with the asymptotics of $I_\gamma(\zz)$. Clearly,
\bea
A(\zz) P(\yy;\zz)\ =\ 
e^{\pi i\frac{N(N-1)}4}\ e^{\pi i {N}}\
\prod_{a=1}^N \ (z_{N+a}-z_a)^{-1/2} (y_a-y_{N+a}) + \dots\ .
\eea
Hence, $C$ is given by formula \Ref{constant}. Theorem \ref{thm main} is proved.

\bigskip

{\it E-mail address}:\ {\bf anv@email.unc.edu}

\end{document}